\documentclass{amsart}[1999/11/29 v2.0]

\usepackage{amssymb,latexsym,amsxtra,amscd,verbatim,ifthen}

\providecommand{\qedhere}{}
\newboolean{ams2.0}
\ifthenelse{\equal{\qedhere}{}}%
	{\setboolean{ams2.0}{false}}%
	{\setboolean{ams2.0}{true}}%

\DeclareMathOperator{\Aut}{Aut}

\DeclareMathOperator{\rank}{rank} 
\DeclareMathOperator{\GKdim}{GKdim}
\DeclareMathOperator{\Supp}{Supp}
\DeclareMathOperator{\Pic}{Pic}

\newcommand{\nequiv}{\not \equiv}
\newcommand{\Num}{\text{Num}}
\newcommand{\ZZ}{{\mathbb Z}}
\newcommand{\RR}{{\mathbb R}}

\newcommand{\NN}{{\mathbb N}}

\newcommand{\LL}{{\mathcal L}}
\newcommand{\F}{{\mathcal F}}
\newcommand{\OO}{{\mathcal O}}

\newcommand{\gs}{{\sigma}}

\numberwithin{equation}{section}

\newtheorem{theorem}[equation]{Theorem}
\newtheorem{cor}[equation]{Corollary}
\newtheorem{lemma}[equation]{Lemma}
\newtheorem{prop}[equation]{Proposition}

\theoremstyle{definition}
\newtheorem{definition}[equation]{Definition}
\newtheorem{example}[equation]{Example}

\theoremstyle{remark}
\newtheorem{remark}[equation]{Remark}
\newtheorem{qu}[equation]{Questions}
\newtheorem{qu1}[equation]{Question}

\begin{document}

\setcounter{page}{1}

\title{Criteria for $\sigma$-ampleness}

\author{Dennis S. Keeler}
     \thanks{ 
     Partially supported by NSF grant DMS-9801148.}  
\address{ Department of Mathematics \\ University of Michigan \\ Ann Arbor, MI 48109-1109   }
\email{dskeeler@umich.edu}
\urladdr{http://www-personal.umich.edu/\~{}dskeeler}
 

\ifthenelse{\boolean{ams2.0}}%
	{}%
	{}

%
%

\begin{abstract}
In the noncommutative geometry of Artin, Van~den~Bergh, and others, the twisted
homogeneous coordinate ring is one of the basic constructions. Such a ring is
defined by a $\sigma$-ample divisor, where $\sigma$ is an automorphism of a
projective scheme $X$. Many open questions regarding $\sigma$-ample divisors have
remained.

We derive a relatively simple necessary and sufficient condition for a divisor on
$X$ to be $\sigma$-ample. As a consequence, we show right and left
$\sigma$-ampleness are equivalent and any associated noncommutative homogeneous
coordinate ring must be noetherian and have finite, integral GK-dimension. We
also characterize which automorphisms $\sigma$ yield a $\sigma$-ample divisor. 
\end{abstract}

\maketitle

\section{Introduction} \label{S:intro}

In the past ten years a study of ``noncommutative projective geometry'' has flourished.
By using and generalizing techniques of commutative projective geometry, one can study
certain noncommutative rings
and obtain results for which no purely algebraic proof is known.

The most basic building block of the theory is the twisted homogeneous coordinate ring.
Let $X$ be a projective scheme over an algebraically closed field $k$
 with $\gs$  a scheme automorphism 
and let $\LL$ be an invertible sheaf on $X$.
In \cite{ATV} a twisted version of the homogeneous coordinate ring $B=B(X,\gs,\LL)$ of $X$
was invented with the grading $B = \oplus B_m$ for
\[
	B_m = H^0 (X, \LL \otimes {\LL}^{\gs}  \otimes \dots \otimes {\LL}^{{\gs}^{m-1}} )
\]
where ${\LL}^{\gs} = {\gs}^*{\LL}$ is the pullback of $\LL$. Multiplication on sections is defined by
$a \cdot b = a \otimes b^{{\gs}^m}$ when $a \in B_m$ and $b \in B_n$. 

Soon after their seminal paper, Artin and Van~ den~ Bergh formalized much of the theory of
these twisted homogeneous coordinate rings in \cite{AV}. In the commutative case, 
the most useful homogeneous
coordinate rings are associated with an ample invertible sheaf. A generalization of ampleness
was therefore needed and defined as follows.

An invertible sheaf $\LL$ 
is called right $\gs$-\emph{ample} if for any coherent sheaf ${\mathcal F}$,
\[ 
	H^q(X, {\mathcal F} \otimes \LL \otimes {\LL}^{\gs} 
 	\otimes \dots \otimes {\LL}^{{\gs}^{m-1}}) = 0
\]
for $q > 0$ and $m \gg 0$.
Similarly, $\LL$ is called left $\gs$-\emph{ample} if for any coherent sheaf ${\mathcal F}$, 
\[
	H^q(X, \LL \otimes {\LL}^{\gs} 
 	\otimes \dots \otimes {\LL}^{{\gs}^{m-1}} \otimes {\mathcal F}^{{\gs}^{m}} ) = 0
\]
for $q > 0$ and $m \gg 0$.
A divisor $D$ is called right (resp. left) $\gs$-ample if ${\mathcal O}_X(D)$ is right (resp. left)
$\gs$-ample.  If $\gs$ is the identity automorphism, then these conditions are the same as saying
$\LL$ is ample.
Artin and Van den Bergh proved that if $\LL$ is right (resp. left) $\gs$-ample, 
then $B$ is a finitely generated
right (resp. left) noetherian $k$-algebra \cite{AV}.
  
Twisted homogeneous coordinate rings have been instrumental in the classification of
 rings,
such as the $3$-dimensional Artin-Schelter regular algebras \cite{ATV,St1,St2} 
and the $4$-dimensional Sklyanin algebras \cite{SS}. 
Artin and Stafford
showed that any connected (i.e. $B_0 = k$) graded domain of GK-dimension $2$ generated by
$B_1$ is the twisted homogeneous coordinate ring (up to a finite dimensional vector space)  
of some projective curve $X$, with automorphism
$\gs$ and (left and right) $\gs$-ample $\LL$ \cite{AS}. 
Therefore any such ring is automatically noetherian!

While the concept of noncommutative schemes has grown to encompass more than just
twisted homogeneous coordinate rings (cf.~ \cite{AZ}), they remain a guide for how such a 
scheme ought to behave. However, fundamental open questions about these coordinate rings 
and $\gs$-ample divisors have
persisted for the past decade. In \cite{AV}, the authors derived a simple criterion for a divisor
to be $\gs$-ample in the case $X$ is a curve, a smooth surface, or certain other special cases.
With this criterion, they showed that $B$ must have finite GK-dimension. In other words,
they showed that $B$ has polynomial growth. They ask

\begin{qu} \cite[Question~ 5.19]{AV} 
\begin{enumerate}
	\item What is the extension of our simple criterion to higher dimensions?
	\item Does the existence of a $\gs$-ample divisor imply that $B$ has polynomial growth?
\end{enumerate}
\end{qu}

The second question was  asked again after \cite[Theorem~ 4.1]{AS}.

One would also like to know if the conditions of right and left $\gs$-ampleness are related
and if $B$ could be right noetherian, but not left noetherian.
One might ask for which (commutative) schemes and automorphisms a $\gs$-ample divisor
even exists and if one can be easily found.

In this paper, all these questions will be settled very satisfactorily. We obtain
\begin{theorem}
The following are true for any projective scheme $X$ over an algebraically closed field.
\begin{enumerate}
	\item Right and left $\gs$-ampleness are equivalent. Thus every associated $B$ is 
		(right and left) noetherian.
	\item A projective scheme $X$ has a $\gs$-ample divisor if and only if the action of
		$\gs$ on numerical equivalence classes of divisors is quasi-unipotent
		(cf.~ $\S$\ref{S:nq} for definitions). In this case,
		every ample divisor is $\gs$-ample.
	\item $\GKdim B$ is an integer if $B = B(X, \gs, \LL)$ and $\LL$ is $\gs$-ample. Here
		$\GKdim B$ is the Gel'fand-Kirillov dimension of $B$ in the sense of \cite{KL}.
\end{enumerate}
\end{theorem}
The first two results are handled in $\S$\ref{S:cor}, while the third is covered in
$\S$\ref{S:GK}.
 
These facts are all consequences of
\begin{theorem}[See Remark~ \ref{rem:main}] \label{th:main}
Let $X$ be a projective scheme with automorphism $\gs$. Let $D$ be a Cartier divisor.
$D$ is (right) $\gs$-ample if and only if $\gs$ is quasi-unipotent and 
\[
D + \gs D + \dots + \gs^{m-1} D
\]
is ample for some $m > 0$.
\end{theorem}
This is the ``simple criterion'' which was already known if
$X$ is a smooth surface  \cite[Theorem~ 1.7]{AV}. We obtain the result mainly by use of 
Kleiman's numerical theory of ampleness \cite{K}.

Besides the results above, we derive other corollaries in~ $\S$\ref{S:cor} and find bounds
for the GK-dimension in~ $\S$\ref{S:GK} via Riemann-Roch theorems.
We also examine what happens in the non-quasi-unipotent case and obtain

\begin{theorem}[See Remark~ \ref{rem:expgrowthdone}] \label{th:expgrowth}
Let $X$ be a projective scheme with automorphism $\gs$. Then the following are 
equivalent:
\begin{enumerate}
\item\label{th:expgrowth1} The automorphism $\gs$ is quasi-unipotent.
\item\label{th:expgrowth2} For all ample divisors $D$, 
	 $B(X, \gs, \OO_X(D))$ has finite GK-dimension.
\item\label{th:expgrowth3} For all ample divisors $D$, $B(X, \gs, \OO_X(D))$ is noetherian.
\end{enumerate}
\end{theorem}

This paper is part of my Ph.D. thesis at the University of Michigan, under the direction
of J.T.~ Stafford.

\section{Reductions} \label{S:reductions}

Throughout this paper, we will work in the case of  a projective scheme $X$ over an algebraically
closed base field of arbitrary characteristic. A variety will
mean a reduced, irreducible scheme. All divisors will be \emph{Cartier divisors}
unless otherwise stated.
For a projective scheme, the group of Cartier divisors, modulo linear equivalence, is naturally
isomorphic to the Picard group of invertible sheaves. Since much of our work will entail
intersection theory, we will work from the divisor point of view. 
Several times we use the facts that the ample divisors form a cone, 
that ampleness depends only
on the numerical equivalence class of a divisor, and 
that ampleness is preserved under an
automorphism. Hence the cone of ample divisors and its closure, the
cone of numerically effective divisors, are invariant under an automorphism.
As a reference for these and related facts
we suggest \cite{K}.

\begin{remark}
The main results of this paper will be proved in terms of
divisors rather
than line bundles. However, the reader should note that,
unravelling the definitions, one has  $\OO_X(\sigma D) \cong
\OO_X(D)^{\sigma^{-1}}$. It is therefore notationally more convenient to
 work with a right $\sigma^{-1}$-ample line bundle
$\LL=\OO_X(D)$, since then $D$ is right $\sigma^{-1}$-ample
if and only if $H^q(X, \F\otimes \OO_X(D+\sigma D+\dots + \sigma^{m-1} D))=0$
for all $q>0$ and $m \gg 0$. Obviously, this will have no effect on the final theorems.
Throughout this paper, we will use the notation 
$\Delta_m = D+\sigma D+\dots + \sigma^{m-1} D$.
\end{remark}

Before deriving our main criterion for $\gs$-ampleness, we must first prove other equivalent
criteria. We will need

\begin{lemma}\label{lem:Fujita} \cite[p.~ 520, Theorem~ 1]{Fuj}
Let $\F$ be a coherent sheaf on a projective scheme $X$ and let $H$ be an ample divisor on $X$.
Then there exists an integer $c_0$ such that for all $c \geq c_0$,
\[
	H^q(X, \F \otimes \OO_X(cH + N)) = 0
\]
for $q > 0$ and any ample divisor $N$.\qed
\end{lemma}

\begin{prop} \label{prop:tfae}
Let $X$ be a projective scheme with $\gs$ an automorphism. Let $D$ be a divisor on $X$
and $\Delta_m = D + \gs D + \dots + \gs^{m-1} D$. Then
the following are equivalent:
\begin{enumerate}
	\item \label{prop:tfae1} For any coherent sheaf $\F$, there exists an $m_0$ such that
		\[
		H^q(X, \F \otimes \OO_X(\Delta_m)) = 0
		\] for  $q > 0$ and $m \geq m_0$.
	\item \label{prop:tfae2} For any coherent sheaf $\F$, there exists an $m_0$ such that
		$\F \otimes \OO_X(\Delta_m)$ is generated by global sections for $m \geq m_0$.
	\item \label{prop:tfae3} For any divisor $H$, there exists an $m_0$ such that
		$\Delta_m - H$ is very ample for $m \geq m_0$.
	\item \label{prop:tfae4} For any divisor $H$, there exists an $m_0$ such that
		$\Delta_m - H$ is ample for $m \geq m_0$.
\end{enumerate}
\end{prop}
The first condition is the original definition of right $\gs^{-1}$-ample.
\begin{proof}
 (\ref{prop:tfae1}) $\Rightarrow$ (\ref{prop:tfae2}) is
\cite[Proposition~ 3.2]{AV}.

(\ref{prop:tfae2}) $\Rightarrow$ (\ref{prop:tfae3}) follows from the fact 
that a very ample divisor plus
a divisor generated by its sections is a very ample divisor. Given any divisor $H$ and a very
ample divisor $H'$, choose $m_0$ such that $\Delta_m - H - H'$ is generated by global sections
for $m \geq m_0$. Then $\Delta_m - H - H' + H' = \Delta_m - H$ is very ample for $m \geq m_0$.

(\ref{prop:tfae3}) $\Rightarrow$ (\ref{prop:tfae4}) is trivial.

(\ref{prop:tfae4}) 
$\Rightarrow$
(\ref{prop:tfae1}).
For any ample divisor $H$ and any $c \geq 0$, one can 
choose $m_0$ so that 
for $m \geq m_0$,  we have $N = \Delta_m - cH$ is an ample divisor. 
Then (\ref{prop:tfae1}) follows immediately from Lemma~ \ref{lem:Fujita}.
\end{proof}

A similar proposition holds for left $\gs^{-1}$-ample divisors, with $\F$ and $H$ replaced by
$\F^{\gs^{-m}}$ and $\gs^{m} H$. One deduces this easily from

\begin{lemma} \label{lem:rightleft} \cite[p.~ 31]{Ste}
A divisor $D$ is right $\gs^{-1}$-ample if and only if $D$ is left $\gs$-ample.
\end{lemma}
\begin{proof}
Let $D$ be right $\gs^{-1}$-ample. Then
for any coherent sheaf $\F$, there exists an
$m_0$ such that
\[
H^q(X, \OO_X(D + \gs D + \dots + \gs^{m-1} D) \otimes \F^{\gs}) = 0 
\] for $q > 0$ 
and $m \geq m_0$. Since cohomology
is preserved under automorphisms, pulling back by $\gs^{m-1}$, we have
\[
H^q( X, \OO_X( D + \gs^{-1} D + \dots + \gs^{-(m-1)} D) \otimes \F^{\gs^{m}}) = 0
\]
for $q > 0$ and $m \geq m_0$. So $D$ is left $\gs$-ample.
\end{proof}

It is often useful to replace $D$ with $\Delta_m$ and $\gs$ with $\gs^m$ to assume
$D$ and $\gs$ have a desired property. Using standard techniques, one can also show

\begin{lemma} \label{lem:replace} \cite[Lemma~ 4.1]{AV}
Let $D$ be a divisor on $X$. Given a positive integer $m$,
 $D$ is right $\gs^{-1}$-ample if and only if $\Delta_m$ is right ${\gs}^{-m}$-ample.\qed
\end{lemma}

\section{The non-quasi-unipotent case} \label{S:nq}

Let $A^1_{\Num}(X)$ be the set of divisors of $X$ modulo numerical
equivalence. That is, for divisors $D$ and $D'$, one has
$D \equiv D' \in A^1_{\Num}(X)$ if and only if $(D.C) = (D'.C)$ for all integral
curves $C \subset X$. We will use this definition implicitly several times,
especially the fact that a non-zero element of $A^1_{\Num}(X)$ has
non-zero intersection with some curve.
Further, one has that $A^1_{\Num}(X)$ is a finitely generated free abelian group 
\cite[p.~ 305, Remark~ 3]{K}.
 Let $P$ be the action of $\gs$ on $A^1_{\Num}(X)$; hence $P \in 
{\mathrm{GL}}(\ZZ^{\ell})$ for some $\ell$.

A matrix is called 
\emph{quasi-unipotent} if all of its eigenvalues
are roots of unity. We call an automorphism $\gs$ quasi-unipotent if $P$ is.
The main goal of this section is to show that a non-quasi-unipotent $\gs$ cannot give a
$\gs$-ample divisor.

First, we must review an useful fact about integer matrices.
\begin{lemma}
Let $P \in {\mathrm{GL}}(\ZZ^\ell)$. Then $P$ is quasi-unipotent if and only if all eigenvalues of
$P$ have absolute value $1$. Thus if $P$ is not quasi-unipotent, then $P$ has an eigenvalue
of absolute value greater than $1$.
\end{lemma}
\begin{proof}
The first claim is \cite[Lemma~ 5.3]{AV}.
For the second claim, 
the property of $P$ not being quasi-unipotent is reduced to saying $P$ has an eigenvalue of
absolute value not $1$. Since $P$ has determinant $\pm 1$, $P$ has an eigenvalue of absolute value
greater than $1$.
\end{proof}

The following lemma shows a relationship between 
the spectral radius $r = \rho(P)$ and the intersection numbers
$(\gs^m D. C)$, where $D$ is an ample divisor.

\begin{lemma} \label{lem:pairgrowth}
Let $P$ be as described above with spectral radius $r = \rho(P)$.
There exists an integral curve $C$ with the following property:
If $D$ is an ample divisor, then there exists  $c > 0$ such that
\[
	(\gs^m D.C) \geq c r^{m} \qquad \text{for all }  m \geq 0.
\]
\end{lemma}
\begin{proof}
Let $\kappa$ be the cone generated by numerically effective divisors
in $A^1_{\Num}(X) \otimes \RR$. In the terminology of
\cite{V}, $\kappa$ is a solid cone since it has a non-empty interior \cite[p.~ 325, Theorem~ 1]{K}.
Since $P$ maps $\kappa$ to $\kappa$, the spectral radius $r$ is an eigenvalue of $P$ and $r$ has an
eigenvector $v \in \kappa$ \cite[Theorem~ 3.1]{V}.

Since $v \in \kappa \setminus \{0\}$, 
there exists a curve $C$ with $(v. C) > 0$. 
Given an ample divisor
$D$, there is a positive $\ell$ so that $\ell D - v$ is in the ample cone \cite[p.~1209]{V}. 
Thus
\[
\ell (\gs^m D. C) = \ell (P^m D. C) > (P^m v. C) = r^m (v. C).
\]
Taking $c = (v. C)/\ell$, we have the lemma.
\end{proof}

Now a graded ring $B = \oplus_{i \geq 0} B_i$ is finitely graded 
if  $\dim B_i < \infty$ for all $i$.
Such a ring $B$ has exponential growth (see \cite{SZ}) if
\begin{equation}
\limsup_{n \to \infty}
 \biggl(\sum_{i \leq n} \dim B_i \biggr)^{\frac{1}{n}} > 1.
\end{equation}
If $B$ has exponential growth, it  
is neither right nor left noetherian  \cite[Theorem~0.1]{SZ}. 
This fact combined with the intersection numbers above allow us to prove

\begin{theorem}\label{th:sampleimpliesquasi}
Let $X$ be a projective scheme with automorphism $\gs$. 
If $X$ has a right $\gs^{-1}$-ample divisor, then $\gs$ is quasi-unipotent.
\end{theorem}
\begin{proof}
Suppose that $D$ is a right $\gs^{-1}$-ample
divisor. Let $\Delta_m = D + {\gs}D + \dots + {\gs}^{m-1}D$. By (\ref{prop:tfae}) and
(\ref{lem:replace}), we may replace $D$ with $\Delta_m$ and $\sigma$ with $\sigma^m$
and assume that $D$ is ample.

Let $P$ be the action of $\gs$ on $A^1_{\Num}(X)$.
Suppose $P$ is non-quasi-unipotent with spectral radius $r > 1$ and choose an integral
curve $C$ as in Lemma~\ref{lem:pairgrowth}. Let ${\mathcal I}$ be the ideal sheaf
defining $C$ in $X$. Since $D$ is right $\gs^{-1}$-ample, the higher cohomologies of
${\mathcal I}(\Delta_m) = {\mathcal I} \otimes \OO_X(\Delta_m)$ 
and $\OO_C(\Delta_m)$ vanish for $m \gg 0$. So one has
an exact sequence
\[
0 \to H^0(X, {\mathcal I}(\Delta_m)) \to H^0(X, \OO_X(\Delta_m)) \to
H^0(C, \OO_C(\Delta_m)) \to 0.
\]

For $m \gg 0$, the Riemann-Roch formula for curves \cite[p.~360, Example~18.3.4]{Fu} gives
\[
\dim H^0(C, \OO_C(\Delta_m)) = (\Delta_m. C) + \text{a constant term}.
\]
Thus using the exact sequence and the previous lemma, there exists $c > 0$ so that
\[
\dim H^0(X, \OO_X(\Delta_m)) > c r^m
\]
for $m \gg 0$.
Thus the associated twisted homogenous coordinate ring has exponential growth
and hence is not (right or left) noetherian \cite[Theorem~ 0.1]{SZ}. So
$D$ cannot be right $\gs^{-1}$-ample, by \cite[Theorem~1.4]{AV}.
\end{proof}

\begin{remark}
One can give a more elementary, computational proof of 
Theorem~\ref{th:sampleimpliesquasi}.
Indeed, examining the Jordan form of $P$ gives an upper bound
on $(\gs^m D. C)$.
 Further, using the full strength of 
\cite[Theorem~ 3.1]{V} and asymptotic estimates, 
one can improve the lower bound of Lemma~\ref{lem:pairgrowth}.
We then have
\begin{equation}\label{eq:betterbounds}
c_1 m^k r^m > (\gs^m D. C) > c m^k r^m
\end{equation}
for $m > 0$, where $k+1$ is the size of the largest Jordan block
associated to $r$.
Then using estimates similar to those in the proof of
\cite[Lemma~ 5.10]{AV}, one can find an ample divisor $H$ such that
\[
(\Delta_m - H. \gs^{m}C) < 0
\]
for all $m \gg 0$. This contradicts the
fourth equivalent condition for
right $\gs^{-1}$-ampleness in Proposition~\ref{prop:tfae}.
\end{remark}

Even when an automorphism $\gs$ is not quasi-unipotent, one can form associated
twisted homogeneous coordinate rings. As might be expected, some of these rings
have exponential growth.

\begin{prop}\label{prop:expgrowth}
Let $X$ be a projective scheme with non-quasi-unipotent automorphism $\gs$.
Let $D$ be an ample divisor. Then there exists an integer $n_0 > 0$ such that
for all $n \geq n_0$, the ring
$B = B(X, \gs, \OO_X(nD))$ has exponential growth and is neither
right nor left noetherian.
\end{prop}
\begin{proof}
Again choose a curve $C$ as in Lemma~\ref{lem:pairgrowth}
with ideal sheaf ${\mathcal I}$. By Lemma~\ref{lem:Fujita}, there exists
$n_0$ such that for all
$n \geq n_0$ and $q > 0$,
\[
H^q(X, {\mathcal I}(nD + N)) = H^q(C, \OO_C(nD + N)) = 0
\]
for any ample divisor $N$. In particular, the above cohomologies vanishes
for $nD + N = nD + \gs(nD) + \dots + \gs^{m-1}(nD)$ where $m > 1$.
Then repeating the last paragraph of 
the proof of Theorem~\ref{th:sampleimpliesquasi} shows that $B$ has
exponential growth.
\end{proof}

When $X$ is a nonsingular surface, \cite[Corollary~5.17]{AV} shows that the above
proposition is true for $n_0 =1$. Their proof makes use of the relatively simple
form of the Riemann-Roch formula and the vanishing of $H^{2}(X, \OO_X(\Delta_m))$
when $\Delta_m$ is the sum of sufficiently many ample divisors. The proof
easily generalizes to the singular surface case, but not to higher dimensions.

\begin{qu1}
Given a non-quasi-unipotent automorphism $\gs$ and ample divisor $D$ on
a scheme $X$, must $B(X, \gs, \OO_X(D))$ have exponential growth?
\end{qu1}

There do exist varieties with non-quasi-unipotent automorphisms. 
If the canonical divisor $K$ is
ample or minus ample, 
then  any automorphism $\gs$ must be quasi-unipotent (cf. (\ref{prop:quasiexamples})). 
So intuitively, one expects to find non-quasi-unipotent automorphisms far away
from this case, i.e., when $K=0$. Further, there are strong existence theorems for
automorphisms of $K3$ surfaces (which do have
$K=0$). Indeed, a $K3$ surface with non-quasi-unipotent automorphism is studied in \cite{W}.

\begin{example}
\emph{There exists a $K3$ surface with automorphism $\gs$ such that $X$ has no $\gs$-ample divisors.}

\begin{proof}
Wehler \cite[Proposition~2.6, Theorem~ 2.9]{W}
constructs a family of $K3$ surfaces whose general member 
$X$ has
\[
	\Pic(X) \cong A^1_{\Num}(X) \cong \ZZ^2, \qquad \Aut(X) \cong \ZZ/2\ZZ \ast \ZZ/2\ZZ.
\]
(That is, $\Aut(X)$ is the free product of two cyclic groups of order 2.)
The ample generators $H_1$ and $H_2$ of $A^1_{\Num}(X)$ have intersection numbers
\[
(H_1^2) = (H_2^2) = 2, \qquad (H_1.H_2)=4.
\]

$\Aut(X)$ has two generators $\gs_1, \gs_2$ whose actions on $A^1_{\Num}(X)$ can be
represented as two quasi-unipotent matrices
\[
	\gs_1 =
		\begin{pmatrix}
			1 & 4 \\
			0 & -1
		\end{pmatrix},
		\qquad
	\gs_2 =
		\begin{pmatrix}
			-1 & 0 \\
			4  & 1
		\end{pmatrix}.
\]
However, the action of $\gs_1 \gs_2$ has eigenvalues $7 \pm 4 \sqrt{3}$. So 
\emph{$X$ has no
$\gs_1 \gs_2$-ample divisor}. Note that by Corollary~ \ref{cor:ample_sigmaample} below,
any ample divisor is $\gs_1$-ample and
$\gs_2$-ample.
\end{proof}
\end{example}

\section{The quasi-unipotent case} \label{S:q}

Now let $\gs$ be a quasi-unipotent automorphism with $P$ its action on $A^1_{\Num}(X)$.
We will have 
several uses for a particular invariant of $\gs$.

\begin{definition}\label{def:J}
Let $k+1$ be the rank of the largest Jordan block of $P$.
We define $J(\gs)=k$.
\end{definition}

Note that  $J(\gs)=J(\gs^m)$ for all $m \in \ZZ \setminus \{0\}$.
 It may be that $k$ is greater than $0$, as seen in 
\cite[Example~ 5.18]{AV}. 
We will see in the next section that $k$ must be even, but this
is not used here.

To prove Theorem~ \ref{th:main}, it remains to
show that (for $\gs$ quasi-unipotent) if $D$ is a divisor
such that $\Delta_m$ is ample for some $m$, then
$D$ is right $\gs^{-1}$-ample. So fix such a $D$. We may
again replace $D$ with $\Delta_n$ and $\gs$ with $\gs^n$ via
(\ref{lem:replace}), so that $D$ is ample and $P$ is
unipotent, that is $P=I+N$, where $N$ is the nilpotent part of $P$.
In this case, $k = J(\gs)$ is the smallest natural number
such that $N^{k+1} = 0$.

 We let $\equiv$ denote numerical equivalence and reserve $=$ for linear equivalence. We then
 have, for all $m \geq 0$,
\begin{align}
\label{eq:sigmaD}	 \sigma^mD &\equiv P^m D = \sum_{i=0}^k \binom{m}{i} N^iD, \\
\label{eq:Delta_m} 
	\Delta_m &\equiv \sum_{i=0}^k \binom{m}{i+1} N^iD.
\end{align}

Once a basis for $A^1_{\Num}(X)$ is chosen, one can treat $N^i D$ as a divisor. Of course,
this representation of $N^i D$ is not canonical. However, since ampleness and intersection 
numbers only depend on numerical equivalence classes, this is not a problem.


\begin{lemma}\label{lem:notzero}
Let $\gs$ be a unipotent automorphism with $P = I + N$ and $k = J(\gs)$. If $D$ is an
ample divisor, then $N^k D \nequiv 0$ in $A^1_{\Num}(X)$.
\end{lemma}
\begin{proof}
Since $N^k \neq 0$, there exists a divisor $E$ and
curve $C$ such that $(N^k E. C) > 0$. Choose $\ell$ so that $\ell D - E$ is ample.
By Equation~\ref{eq:sigmaD}, the intersection numbers $(\gs^m(\ell D - E).C)$ are given
by a polynomial in $m$ with leading coefficient $(\ell N^k D - N^kE . C) / k!$. Since this
polynomial must have positive values for all $m$, we must have $N^kD \nequiv 0$.
\end{proof}

We now turn towards proving that for any divisor $H$, there exists $m_0$ such that 
$\Delta_{m_0} - H$ is ample,
when $\gs$ is unipotent and $D$ is ample. Then since $D$ is ample,
$\Delta_m - H$ is ample for $m \geq m_0$. For certain $H$, this is true even if
 $\gs$ is not quasi-unipotent.
 
\begin{lemma}
Let $X$ be a projective scheme with automorphism $\gs$ (not necessarily
quasi-unipotent). Let $D$ be an ample divisor and $H$ a divisor whose numerical equivalence class
is fixed by $\gs$. Then there exists an $m$ such that $\Delta_m - H$ is ample.
\end{lemma}
\begin{proof}
Choose $m$ such that $D' = mD - H$ is ample. Let 
\[
\Delta_j' = D' + \gs D' + \dots + \gs^{j-1} D'.
\]
Then
$\Delta_m' \equiv m \Delta_m - mH$
is ample and thus $\Delta_m - H$ is ample.
\end{proof}

\begin{prop}\label{prop:unipotentgivesample}
Let $X$ be a projective scheme with unipotent automorphism $\gs$. Let $D$ be an ample
divisor and $H$ any divisor. Then there exists an $m_0$ such that $\Delta_{m_0} - H$ is ample.
Hence $\Delta_m - H$ is ample for $m \geq m_0$.
\end{prop}
\begin{proof}
Let $W \subset A^1_{\Num}(X) \otimes \RR$ be the span of $D, ND, \dots, N^k D$.
$W$ is a $k+1$-dimensional vector space by Lemma~ \ref{lem:notzero}. 
By Equation~ \ref{eq:sigmaD},
 it contains the real cone $\kappa$ generated by
$S = \{ \gs^i D \vert i \in \NN \}$. Using a lemma of Caratheodory \cite[p.~ 45, Lemma~ 1]{Ha2}, any
element of $\kappa$ can be written as a linear combination of $k+1$
elements of $S$ with non-negative real coefficients. Thus for all $m \in \NN$,
\[
	 \Delta_m \equiv \sum^{k}_{i=0} f_i(m) \gs^{g_i(m)} D
\]
where $f_i\colon \NN \to \RR_{\geq 0}$ and $g_i\colon \NN \to \NN$. 
Expanding the $\gs^{g_i(m)}D$ above
and comparing the coefficient of $D$
 with Equation~ \ref{eq:Delta_m}, one finds that
\[
	\sum^k_{i=0} f_i(m) = m.
\]
Since $f_i(m) \geq 0$, for each $m$, 
there must be some $j$ such that $f_j(m) \geq m/(k+1)$.

Now choose $l$ such that $lD - H$ is ample and $m_0$ such that, $m_0/(k+1) \geq l$. Then
\[
	f_j(m_0) \gs^{g_j(m_0)} D - \gs^{g_j(m_0)}H
\]
 is in the ample cone for the given $j$. Set $g=g_j(m_0)$.
 The other $f_i(m_0)$ are non-negative. 
Then $\Delta_{m_0} - \gs^{g} H$ is in the ample cone as it is a sum of elements in the ample cone. 
 Since
 it is a divisor, it is ample \cite[p.~ 324, Remark~ 3]{K}.
 
We now prove the lemma by induction on $q$, the smallest positive integer such that
$N^q H \equiv 0$. Since $N$ is nilpotent, there is such a $q$ for any $H$. The case $q=1$ is
handled by the previous lemma.

Now as $\gs \equiv I + N$, we know $\gs^{-m_0}(\gs^{g} H - H)$ is killed by $N^{q-1}$. 
So there is an $m_1$ so that
\[
	Y= \Delta_{m_1} + \gs^{-m_0}(\gs^{g} H - H)
\]
 is ample. Then as $\gs$ fixes the ample cone
\[
	\Delta_{m_0} - \gs^{g} H + \gs^{m_0} Y = \Delta_{m_0 + m_1} - H
\]
is ample.
\end{proof}

We now immediately have by Propositions~ \ref{prop:tfae}, \ref{prop:unipotentgivesample},
 and Theorem~ \ref{th:sampleimpliesquasi}:

\begin{theorem} \label{th:sigmaample}
Let $X$ be a projective scheme with an automorphism $\gs$.
A divisor $D$ is right $\gs^{-1}$-ample if and only if
$\gs$ is quasi-unipotent and
$D + \gs D + \dots + \gs^{m-1} D$ is ample for some $m$.\qed
\end{theorem}

\section{Corollaries}\label{S:cor}

The characterization of (right) $\gs^{-1}$-ampleness has many strong corollaries which are now easy
to prove, but were only conjectured before.

\begin{cor}
Right $\gs$-ample and left $\gs$-ample are equivalent conditions. Further, $\gs$-ampleness and
$\gs^{-1}$-ampleness are equivalent.
\end{cor}

\begin{proof}
Let $D$ be right $\gs^{-1}$-ample. By Theorem~ \ref{th:sigmaample}, $\gs$ is quasi-unipotent and
$\Delta_m$ is ample for some $m$. 
Then $\gs^{-1}$ is quasi-unipotent and 
\[
\gs^{-(m-1)} \Delta_m = D + \gs^{-1} D + \dots + \gs^{-(m-1)} D
\]
is ample. Applying the theorem again,
we have that $D$ is right $\gs$-ample. Thus $D$ is left $\gs^{-1}$-ample by
Lemma~ \ref{lem:rightleft}. The same lemma gives the second statement of the corollary.
\end{proof}

\begin{remark} \label{rem:main}
Combined with (\ref{th:sigmaample}), this proves Theorem~ \ref{th:main}
and so we may refer to a divisor as being simply ``$\gs$-ample.''
\end{remark}

 In \cite{AV}, left
$\gs$-ampleness was shown 
to imply the associated twisted homogeneous coordinate
ring is left noetherian.
However, as noted in the footnote of \cite[p.~ 258]{AS}, the paper says, but does not prove,
that $B$ is noetherian. 
This actually is the case.
\begin{cor}\label{cor:noetherian}
Let $B=B(X, \gs, \OO_X(D))$ be the twisted homogeneous coordinate ring associated to a
$\gs$-ample divisor $D$. Then $B$ is a (left and right) noetherian ring, 
finitely generated over the base field.\qed
\end{cor}

Analysis of the GK-dimension of $B$ will be saved for the next section.

From the definition of $\gs$-ample, it is not obvious when $\gs$-ample divisors even exist.
Theorem~ \ref{th:sigmaample} makes the question much easier.

\begin{cor}\label{cor:ample_sigmaample}
A projective scheme $X$ has a $\gs$-ample divisor if and only if $\gs$ is quasi-unipotent. In
particular, every ample divisor is a $\gs$-ample divisor if $\gs$ is quasi-unipotent.\qed
\end{cor}

Thus, it is important to know when an automorphism $\gs$ is 
quasi-unipotent. 
From the bounds in Equation~\ref{eq:betterbounds}, we obtain

\begin{prop}
Let $D$ be an ample divisor. Then $\gs$ is quasi-unipotent if and only if
for all curves $C$, the intersection numbers $(\gs^m D. C)$ are bounded
by a polynomial for positive $m$.\qed
\end{prop}

\begin{prop}\label{prop:quasiexamples}
Let $X$ be a projective scheme such that
\begin{enumerate}
\item $X$ has a canonical divisor $K$ which is an ample
or minus-ample divisor, or
\item the Picard number of $X$, i.e., the rank of $A^1_{\Num}(X)$, is $1$.
\end{enumerate}
Then any automorphism $\gs$ of $X$ is quasi-unipotent.
Indeed, some power of $\gs$ is numerically equivalent to the identity.
\end{prop}
\begin{proof}
In the first case, for $K$ to be ample or minus-ample, it must be a
Cartier divisor. Thus the intersection numbers $(\gs^m K. C)$ are defined,
where $C$ is a curve. Since $K$ must be fixed by $\gs$, some power
of $\gs$ must be numerically equivalent to the identity by
Equation~\ref{eq:betterbounds}. In the second case, the
action of $\gs$ itself must be numerically equivalent to the identity.
\end{proof}

Thus for many important projective varieties, such as curves, projective $n$-space,
Grassmann varieties \cite[p.~271]{Fu}, and Fano varieties
 \cite[p.~240, Definition~1.1]{KolRat}, one automatically
has that any automorphism must be quasi-unipotent.

Returning to corollaries of
Theorem~\ref{th:sigmaample}, 
we see that building new $\gs$-ample divisors from old ones is also possible.

\begin{cor}
Let $D$ be a $\gs$-ample divisor and let $D'$ be a divisor with one of the following properties:
\begin{enumerate}
	\item $\gs$-ample,
	\item generated by global sections, or
	\item numerically effective.
\end{enumerate}
Then $D + D'$ is $\gs$-ample.
\end{cor}
\begin{proof}
Take $m$ such that $\Delta_m$ is ample and $\Delta_m' = D' + \dots + \gs^{m-1} D'$ 
is respectively ample, generated by global sections,
or numerically effective. 
Then $\Delta_m + \Delta_m'$ is ample and we again apply the main theorem.
\end{proof}

The following could be shown directly from the definition, 
but also using a similar method to the proof
above, one can see
\begin{cor}
Let $\gs$ and $\tau$ be automorphisms. Then $D$ is $\gs$-ample if and only if $\tau D$ is
$\tau \gs \tau^{-1}$-ample.\qed
\end{cor}
Note that $\tau$ need not be quasi-unipotent.

Finally, as in the case of ampleness, $\gs$-ampleness is a numerical condition.
\begin{cor}
Let $D, D'$ be numerically equivalent divisors and $\gs, \gs'$
 be numerically equivalent automorphisms (i.e., their actions on $A^1_{\Num}(X)$ are equal).
Then  $D$ is $\gs$-ample if and only if $D'$ is $\gs'$-ample.
\end{cor}
\begin{proof}
As $\Delta_m \equiv D' + (\gs') D' + \dots + (\gs')^{m-1} D'$ and ampleness depends only on 
the numerical equivalence class of a divisor, the corollary follows from our main theorem.
\end{proof}

\section{GK-dimension of $B$} \label{S:GK}

As mentioned, our main goal of this section is to prove

\begin{theorem} \label{th:GKdim}
Let $B = B(X, \gs, \LL)$ for some projective scheme $X$ and $\gs$-ample invertible 
sheaf $\LL$.
\begin{enumerate}
\item \label{th:GKdim0} $\GKdim B$ is an integer. Hence $B$ is of polynomial growth.
In addition, $\GKdim B$ is independent
of the $\gs$-ample $\LL$ chosen.
\item \label{th:GKdim1} If $\gs^m \equiv I$  for some $m$, then
$\GKdim B = \dim X + 1$.
\item \label{th:GKdim2} If $X$ is an equidimensional scheme,
\[
k + \dim X + 1 \leq \GKdim B \leq k(\dim X - 1) + \dim X + 1
\]
 where $k = J(\gs)$ (cf. (\ref{def:J})) is an even natural
number depending only on $\gs$.
\end{enumerate}
\end{theorem}

\begin{remark}\label{rem:expgrowthdone}
We now have all the necessary pieces of Theorem~\ref{th:expgrowth}.
In that result, 
(\ref{th:expgrowth1}) $\implies$ (\ref{th:expgrowth2}) is the theorem above.
(\ref{th:expgrowth1}) $\implies$ (\ref{th:expgrowth3}) is from
Corollaries~\ref{cor:noetherian} and \ref{cor:ample_sigmaample}.
And finally (\ref{th:expgrowth2}) $\implies$ (\ref{th:expgrowth1}) and
(\ref{th:expgrowth3}) $\implies$ (\ref{th:expgrowth1}) both follow from
Proposition~\ref{prop:expgrowth}.
\end{remark}

Theorem~\ref{th:GKdim} generalizes \cite[Proposition~ 1.5, Theorem~ 1.7]{AV}.
The authors of \cite{AV} further show that if $X$ is a smooth surface, then $k = 0, 2$ 
and thus the only
possible GK-dimensions are $3$ and $5$.
The proof that $k \leq 2$ in the surface case
 uses the Hodge Index Theorem and thus far we have been
unable to find a similar bound in higher dimensions.
 Note that if 
$X$ is a curve or $X = {\mathbb P}^n$, then $\rank A^1_{\Num}(X) = 1$ and hence 
by Proposition~\ref{prop:quasiexamples}, some power of $\gs$ is numerically
equivalent to the identity (in fact, $P = I$).
So the theorem implies that $\GKdim B = \dim X + 1$.

In studying the GK-dimension of $B = B(X, \gs, \OO_X(D))$ with
$D$ $\gs$-ample, \cite[p.~ 263]{AV} proves that
\begin{equation} \label{eq:Veronese}
\GKdim B(X, \gs, \OO_X(D)) = \GKdim B(X, \gs^m, \OO_X(\Delta_m))
\end{equation}
 for any positive $m$. 
Therefore, we may again assume $P$ is unipotent,
$D$ is  ample, and $H^q(X, \OO_X(\Delta_m)) = 0$ for $q > 0$ and \emph{all} $m>0$. Then
\[
\dim B_m = \dim H^0(X, \OO_X(\Delta_m)) = \chi(\OO_X(\Delta_m))
\]
where $\chi$ is the Euler
characteristic on $X$. We will soon see this is a polynomial in $m$ with positive leading
coefficient. Again replacing $B$ with an appropriate Veronese subring, we may assume
$\dim B_m \leq \dim B_{m+1}$ for all $m \geq 0$. Then the proof of 
\cite[Lemma~ 1.6]{AS} shows that
\begin{equation} \label{eq:GK}
\GKdim B = \deg (\dim B_m) + 1 = \deg (\chi(\OO_X(\Delta_m))) + 1.
\end{equation}

Thus far, we have only used the intersection numbers $(D.C)$, where $D$ is a divisor and
$C$ is a curve.
In studying the growth of $\Delta_m$ in terms of $m$, we will need to examine
the intersection of divisors on higher dimensional subvarieties.
More precisely, for an $n$-dimensional variety $V$, we use the
symmetric $n$-linear form 
\[
(D_1. \dots . D_n)_V = (\OO_X(D_1). \dots . \OO_X(D_n) . \OO_V)
\]
defined in \cite[p.~296, Definition~1]{K}.

Recall 
that a polynomial with rational coefficients, integer valued on integers,
 is called a \emph{numerical polynomial}. We prove

\begin{lemma} \label{lem:polygrowth}
Let $X$ be a projective scheme with unipotent
  automorphism $\gs$ and ample divisors
 $D$ and $D'$ with
 $\Delta_m' = D' + \dots + \gs^{m-1} D'$.  
Further let $V$ be a closed subvariety of $X$ of dimension $n$. 
Then for $0 \leq i \leq n$,
\begin{enumerate}
\item \label{lem:polygrowth1} $(D^i.\Delta_m^{n-i})_V$ is a numerical polynomial
in $m$ with positive leading
coefficient.
\item \label{lem:polygrowth2} $\deg (D^i.\Delta_m^{n-i})_V = \deg ((D')^i.(\Delta_m')^{n-i})_V$.
\item \label{lem:polygrowth3} 
$\deg (D^{i-1}.\Delta_m^{n-i})_W \leq \deg (D^i.\Delta_m^{n-i})_V$
where  $W \subset V$ is a closed subvariety with $\dim W = \dim V - 1$.
\item \label{lem:polygrowth4} $\deg (D^i.\Delta_m^{n-i})_V < \deg (D^{i-1}.\Delta_m^{n-i+1})_V$.
\item \label{lem:polygrowth5} $\deg (\Delta_m^j)_W < \deg (\Delta_m^n)_V$
where $W \subset V$ is a closed subvariety and $\dim W = j < n$.
\end{enumerate}
\end{lemma}

\begin{proof}
Since $\gs$ is unipotent and intersection numbers only depend on numerical equivalence
classes, we may replace $\Delta_m$ by the divisor on the right hand side of 
Equation~\ref{eq:Delta_m}. As noted below that equation, it is not a problem
to treat the $N^iD$ as divisors.
Since the intersection form is multilinear and integer valued
on divisors,
$(D^i.\Delta_m^{n-i})_V$ must be a numerical polynomial.
By the Nakai criterion for ampleness \cite[p.~ 30, Theorem~ 5.1]{Ha2}, 
the function is positive for all positive $m$
(since $\Delta_m$ is ample) and hence 
has a positive leading coefficient.
Thus part~\ref{lem:polygrowth1} is proven.

Now for some fixed $\ell$, we know that $\ell D' - D$ is ample. Hence
\[ 
\ell(D'.D^{i-1}.\Delta_m^{n-i})_V - (D^i.\Delta_m^{n-i})_V =
 (\ell D' - D.D^{i-1}.\Delta_m^{n-i})_V > 0
\]
 for all $m > 0$. Thus
\[
	\deg (D'.D^{i-1}.\Delta_m^{n-i})_V \geq \deg (D^i.\Delta_m^{n-i})_V
\]
and by symmetry the two degrees are equal.
 We can continue this argument, replacing each $D$ with
 $D'$, so $\deg (D^i.\Delta_m^{n-i}) = \deg ((D')^i.\Delta_m^{n-i})$. 
 By also noting that
\[
\ell \Delta_m' - \Delta_m = (\ell D' - D) + \dots + \gs^{m-1} (\ell D' - D)
\]
is ample, one can similarly replace each $\Delta_m$ with $\Delta_m'$. Thus
the second claim is proven.

Now let $W \subset V$ be a closed subvariety with $\dim W = \dim V - 1$. One has
\[
(D_{i-1}.\Delta_m^{n-i})_W = (D^{i-1}.W.\Delta_m^{n-i})_V
\]
by \cite[p.~298, Proposition~5]{K}.
We claim that for some fixed $\ell$, the intersection number of $\ell D - W$ with any
collection of $n-1$ ample divisors 
is positive.
This is well-known if $V$ is normal so $W$ is a Weil divisor,
so for some $\ell$, the Weil divisor $\ell D - W$ is effective \cite[p.~ 282]{R}.
The general case can be seen by pulling back to the normalization of $V$. Since
normalization is a finite, birational morphism, ampleness \cite[p.~25, Proposition~4.4]{Ha2} and
intersection numbers \cite[p.~299, Proposition~6]{K} 
are both preserved under pull-back. Thus the claim is proven.
An argument
similar to the proof of part~\ref{lem:polygrowth2} proves the third claim of the lemma.

For part~\ref{lem:polygrowth4}, Equation~\ref{eq:Delta_m} shows that the leading coefficient of 
$(D^{i-1}.D'.\Delta_m^{n-i})_V$ is a sum
of terms
\[
a_\alpha (D^{i-1}.D'.N^{\alpha_1}D. \dots . N^{\alpha_{n-i}}D)_V
\]
where $a_\alpha ((k+1)!)^n$ is an integer. So any leading coefficient
times $((k+1)!)^n$ is a positive integer.
Thus given any set of ample divisors $\{ D' \}$, there is a $D'$ in that set 
such that $(D^{i-1}.D'.\Delta_m^{n-i})_V$ has the smallest leading coefficient.

Now let 
$j$ be a natural number such that $(D^{i-1}.\gs^jD.\Delta_m^{n-i})_V$ has
the smallest leading coefficient of all $(D^{i-1}.\gs^lD.\Delta_m^{n-i})_V$. 
Then
for any $l > 0$,
\[
\frac{(D^{i-1}.\gs^lD.\Delta_m^{n-i})_V}
	{(D^{i-1}.\gs^jD.\Delta_m^{n-i})_V}
\]
is a rational function with limit, as $m \to  \infty$, greater than or equal to $1$. 
So given any
natural number $M$, 
\[
	 \lim_{m \to \infty} \frac{(D^{i-1}.\Delta_m.\Delta_m^{n-i})_V}
	 					{(D^{i-1}.\sigma^jD.\Delta_m^{n-i})_V} \geq M.
\]

Since this is true for any $M$, the limit must be $+ \infty$. So
\[
	\deg (D^{i-1}.\Delta_m.\Delta_m^{n-i})_V > 
	\deg (D^{i-1}.\gs^jD.\Delta_m^{n-i})_V.
\]
Examining the proof of part~\ref{lem:polygrowth2}, we see the right hand side equals 
$\deg (D^i.\Delta_m^{n-i})_V$, proving part~\ref{lem:polygrowth4}.

Finally, for part~\ref{lem:polygrowth5}, find a chain of subvarieties
$W = V_0 \subsetneq \dots \subsetneq V_{n-j} = V$. Then
part~\ref{lem:polygrowth3} combined with part~\ref{lem:polygrowth4} proves
the claim for each part of the chain.
\end{proof}

By a version of the Riemann-Roch Theorem
for an $n$-dimensional complete scheme $X$ and 
coherent sheaf $\F$ \cite[p.~ 361, Example~18.3.6]{Fu}:
\begin{equation} \label{eq:RR}
\chi(\F(\Delta_m))
 = \sum_{j=0}^n \frac{1}{j!} \int_X (\Delta_m^j) \cap \tau_{X,j}(\F).
\end{equation}
The $\tau_{X,j}(\F)$ is 
a $j$-cycle, a linear combination of $j$-dimensional closed subvarieties
of $\Supp \F$. 
In other words,
\begin{equation} \label{eq:tau}
\tau_{X,j}(\F) = \sum_V a_V [V]
\end{equation}
where $V$ is a subvariety of $X$, $[\thickspace \thickspace ]$ 
denotes rational equivalence, and $a_V$ is a
rational number. The terms of (\ref{eq:RR}), for $\F = \OO_X$ can then be interpreted as
\[
 \int_X (\Delta_m^j) \cap \tau_{X,j}(\OO_X) = \sum_V a_V (\Delta_m^j)_V.
\]

If $X_i$ is  an irreducible component of $X$ of dimension $j$,
then $[X_i] = n [(X_i)_{\text{red}}]$ is a term in 
$\tau_{X,j}(\OO_X)$, where $n$ is the degree of the natural map $(X_i)_{\text{red}} \to X_i$.
To see this, first note that $(\Delta_m^{\dim X_i})_{X_i}/(\dim X_i)!$ must be 
the $\dim X_i$ term of $\chi(\OO_{X_i}(\Delta_m))$ \cite[ibid.]{Fu}. 
Also $a_{(X_i)_{\text{red}}} = n$ by  \cite[p.~298, Corollary~ 2]{K}.
The short exact sequence
\[
0 \to {\mathcal I_i} \to \OO_X \to \OO_{X_i} \to 0
\]
gives $\chi (\OO_X(\Delta_m)) = \chi(\OO_{X_i}(\Delta_m)) + \chi ({\mathcal I_i}(\Delta_m))$.
The support of ${\mathcal I_i}$ does not contain $X_i$ and an
irreducible component is  rationally equivalent only to itself \cite[p.~11, Example~ 1.3.2]{Fu}.
 So there is no $[X_i]$ term in 
$\chi ({\mathcal I_i}(\Delta_m))$ which could cancel out the $[X_i]$ term in the first summand.

\begin{lemma} \label{lem:irreduciblecomponent}
Let $X$ be a projective scheme with unipotent automorphism $\gs$ and
irreducible components $X_i$. 
Then
\[
\deg \chi(\OO_X(\Delta_m)) = \max_{X_i} \deg (\Delta_m^{\dim X_i})_{X_i}.
\]
\end{lemma}

\begin{proof}
If the left hand side is larger than the right hand side, then
by the discussion before the lemma, there is
a subvariety $V$ with
\[
\deg \chi(\OO_X(\Delta_m)) = \deg (\Delta_m^{\dim V})_V > \deg (\Delta_m^{\dim X_j})_{X_j},
\]
where $X_j$ is an irreducible component \emph{properly} containing $V$.
This cannot happen by Lemma~\ref{lem:polygrowth}.\ref{lem:polygrowth5}.

On the other hand, if the right hand side is larger, then there exists
a subvariety $V$ with $a_V < 0$ in the notation of Equation~\ref{eq:tau} and
\[
\deg (\Delta_m^{\dim V})_V = \max_{X_i} \deg (\Delta_m^{\dim X_i})_{X_i}.
\]
The earlier discussion shows that $a_{X_i} > 0$ for each $i$. Hence $V$ is properly
contained in some irreducible component.
But again this cannot happen by (\ref{lem:polygrowth}.\ref{lem:polygrowth5}).
\end{proof}

\begin{lemma} \label{lem:leadingterm}
Let $X$ be a projective scheme with unipotent automorphism $\gs$. 
Let $V \subset X$ be a
closed subscheme which does not contain (the reduction of) an irreducible component of $X$. Then
$\deg \chi(\OO_V(\Delta_m)) < \deg \chi(\OO_X(\Delta_m))$.
\end{lemma}
\begin{proof}
By Lemma~\ref{lem:irreduciblecomponent} we may pick an irreducible component $V_0$ of $V$ with
\[
\deg \chi(\OO_V(\Delta_m)) = \deg (\Delta_m^{\dim V_0})_{V_0}.
\]
Then $X$ has an irreducible component $X_0$ with $V_0 \subsetneq X_0$.
Combining Lemmata~\ref{lem:polygrowth}.\ref{lem:polygrowth5} and 
\ref{lem:irreduciblecomponent}, the claim is proven.
\end{proof}

\begin{prop} \label{prop:numericalpolynomial}
Let $X$ be a projective scheme with unipotent automorphism $\gs$ and ample divisor $D$. 
 Then $\chi(\OO_X(\Delta_m))$ is a numerical
polynomial in $m$. The degree of this polynomial is independent of the ample divisor $D$ chosen. 
Further, if $\gs$ is numerically equivalent to the identity, this polynomial
has degree $\dim X$.
\end{prop}
\begin{proof}
The first claim is obvious since the intersection numbers in Equation~\ref{eq:RR}
are numerical polynomials, as noted in Lemma~\ref{lem:polygrowth}. The independence
of the degree comes from (\ref{lem:polygrowth}.\ref{lem:polygrowth2}).

If $\gs$ is numerically equivalent to the identity, then $k=0$. So $\chi(\OO_X(\Delta_m)) =
\chi(\OO_X(mD))$ has degree $\dim X$. 
\end{proof}

Combined with Equations~ \ref{eq:Veronese} and \ref{eq:GK}, this proposition implies 
the first two parts of Theorem~ \ref{th:GKdim}.

Considering Lemma~\ref{lem:irreduciblecomponent} 
and Equation~\ref{eq:GK}, 
we immediately have
\begin{prop} \label{prop:GKintersection}
Let $X$ be a scheme with unipotent automorphism $\gs$,
 ample divisor $D$, 
and irreducible components $X_i$. 
Let
$B = B(X, \gs, \LL)$. Then
\[
 \GKdim B - 1 
 = \deg \chi(\OO_X(\Delta_m)) = 
 \max_{X_i} \deg (\Delta_m^{\dim X_i})_{X_i}.
\]
In particular, if $X$ is equidimensional, then
\[
\GKdim B - 1 = \deg (\Delta_m^{\dim X})_{X}.
\]
\end{prop}

Note that by replacing $\gs$ by a power, we may assume $\gs$ fixes each irreducible component. 
That is,
$\gs$ is an automorphism of each component. Thus the soon to be proven bounds of 
Theorem~\ref{th:GKdim} for equidimensional schemes can be used to find bounds in the general
case.

\begin{lemma} \label{lem:even}
Let $\gs$ be a unipotent automorphism with numerical action $P = I+ N$, with
$k = J(\gs)$ (cf.~(\ref{def:J})). Then 
$k$ is even and
$\deg \chi(\OO_X(\Delta_m)) \geq k + \dim X$.
\end{lemma}
\begin{proof}
Given an ample divisor $D$, 
one has $N^kD \neq 0$ by Lemma~ \ref{lem:notzero}. So there exists a curve $C$ such that
$(N^kD.C) \neq 0$.
 Since $(\gs^m D.C) > 0$ for all $m \in \ZZ$ and in particular for $m>0$, 
$(N^kD.C) > 0$. However, if $k$ is odd, (\ref{eq:Delta_m}) implies that 
the leading term of $(\gs^{-m}D.C)$ is
 $-\binom{m}{k}(N^kD.C)$ where $m > 0$. Then $(\gs^{-m}D.C) < 0$ for large $m$, which cannot
 occur.
 
For the lower bound, note
 $\deg \chi(\OO_C(\Delta_m)) = \deg (\Delta_m. C) = k+1$. Constructing a chain
 of subvarieties between $C$ and $X$, Lemma~\ref{lem:leadingterm} shows
that $\deg \chi(\OO_X(\Delta_m)) \geq \dim X + k$.
\end{proof}

\begin{lemma} \label{lem:GKupperbound}
Let $n= \dim X$. Then $(\Delta_m^n)_X$ has degree at most $k(n-1) + n$.
\end{lemma}
\begin{proof}
If $k = 0$ the lemma is trivial. So assume that $k>0$. Let $P = I + N$.

Expanding $(\Delta_m^n)$ gives terms of the form
\[
	 f(m)(N^{i_1}D.N^{i_2}D. \dots .N^{i_n}D) 
\]
where $i_1 \leq i_2 \leq \dots \leq i_n$ and $\deg_m f = n + \sum i_j$. We will show that
if $\sum i_j > k(n-1)$ then $(N^{i_1}D.N^{i_2}D. \dots .N^{i_n}D) = 0$.

Order $(i_1, \dots, i_n)$ in the following way: $(i_1, \dots, i_n) > (i_1', \dots, i_n')$ if
the right-most non-zero entry of  $(i_1, \dots, i_n) - (i_1', \dots, i_n')$
 is positive. We proceed by descending induction on this ordering.
 
The largest $n$-tuple in this ordering is $(k, k, \dots, k)$. Since $k > 0$, $N^{k-1}D$ exists
(taking $N^0 = I$) so
\begin{align*}
	 (N^{k-1}D. (N^kD)^{n-1}) &= (P N^{k-1}D. (P N^kD)^{n-1}) \\
	 &= (N^{k-1}D. (N^kD)^{n-1}) + ((N^kD)^n)
\end{align*}
and hence $((N^kD)^n) = 0$.

Now suppose $(i_1, \dots, i_n)$ is such that $\sum i_j > k(n-1)$ and we have proven our claim
for all larger $(i_1', \dots, i_n')$. Since  $\sum i_j > k(n-1)$, we have $i_1 > 0$ so examine
\[
 (N^{i_1-1}D. N^{i_2}D. \dots . N^{i_n}D) = (PN^{i_1-1}D. PN^{i_2}D. \dots . PN^{i_n}D).
\]

A typical term in the right-hand side is of the form
\[
 (N^{i_1-1 + \delta_1}D. N^{i_2 + \delta_2}D. \dots . N^{i_n + \delta_n}D)
\]
where $\delta_j = 0, 1$. The terms with $\delta_j = 1$ where $j> 1$ are all higher in the ordering
than $(i_1, \dots, i_n)$ and hence are zero. This only leaves
\begin{multline*}
(N^{i_1-1}D. N^{i_2}D. \dots . N^{i_n}D) = \\
	(N^{i_1-1}D. N^{i_2}D. \dots . N^{i_n}D) + 
	(N^{i_1}D. N^{i_2}D. \dots . N^{i_n}D)
\end{multline*}
and so $(N^{i_1}D. N^{i_2}D. \dots . N^{i_n}D) = 0$.
\end{proof}

Using Equation~ \ref{eq:Veronese} and Proposition~\ref{prop:GKintersection}, 
these lemmata complete the proof of Theorem~\ref{th:GKdim}.


\begin{example}
Let $X$ be a  $3$-fold and $\gs$ an automorphism with $k = J(\gs) = 2$. 
Then for any $\gs$-ample $D$, 
$\GKdim B = k + \dim X + 1 = 6$. We will not give the proof of this example, since
it is  not particularly illuminating. However, as one might expect, it
consists of expanding $(\Delta_m^3)$ and showing that most of the terms are zero.
\end{example}

\noindent \textbf{Acknowledgements.}
Thanks to J.T.~ Stafford, K.E.~ Smith, R.~ Lazarsfeld, 
     and J.~ Howald for invaluable discussions. Thanks also to the
     referee for finding a mistake in the original manuscript and
     for other helpful suggestions.


\begin{thebibliography}{ATV}

\bibitem[AS]{AS}
M.~Artin and J.~T.~ Stafford, \emph{Noncommutative graded domains with quadratic
  growth}, Invent. Math. \textbf{122} (1995), no.~2, 231--276.

\bibitem[ATV]{ATV}
M.~Artin, J.~Tate, and M.~Van~den Bergh, \emph{Some algebras associated to
  automorphisms of elliptic curves}, The Grothendieck Festschrift, Vol.\ I,
  Birkh\"auser Boston, Boston, MA, 1990, pp.~33--85.

\bibitem[AV]{AV}
M.~Artin and M.~Van~den Bergh, \emph{Twisted homogeneous coordinate rings}, J.
  Algebra \textbf{133} (1990), no.~2, 249--271.

\bibitem[AZ]{AZ}
M.~Artin and J.~J. Zhang, \emph{Noncommutative projective schemes}, Adv. Math.
  \textbf{109} (1994), no.~2, 228--287.

\bibitem[Fj]{Fuj}
Takao Fujita, \emph{Vanishing theorems for semipositive line bundles},
  Algebraic geometry (Tokyo-Kyoto, 1982), Springer, Berlin, 1983, pp.~519--528.

\bibitem[Fl]{Fu}
William Fulton, \emph{Intersection theory}, second ed., Springer-Verlag,
  Berlin, 1998.
  

\bibitem[H]{Ha2}
Robin Hartshorne, \emph{Ample subvarieties of algebraic varieties},
  Springer-Verlag, Berlin, 1970, Notes written in collaboration with C. Musili.
  Lecture Notes in Mathematics, Vol. 156.

\bibitem[K]{K}
Steven~L. Kleiman, \emph{Toward a numerical theory of ampleness}, Ann. of Math.
  (2) \textbf{84} (1966), 293--344.
  
\bibitem[Ko]{KolRat}
J{\'a}nos Koll{\'a}r, \emph{Rational curves on algebraic varieties},
  Springer-Verlag, Berlin, 1996.
  
\bibitem[KL]{KL}
G{\"u}nter~R. Krause and Thomas~H. Lenagan, \emph{Growth of algebras and
  {G}elfand-{K}irillov dimension}, revised ed., American Mathematical Society,
  Providence, RI, 2000.

\bibitem[R]{R}
Miles Reid, \emph{Canonical $3$-folds}, Journ\'ees de G\'eom\'etrie Alg\'ebrique
  d'Angers, Juillet 1979/Algebraic Geometry, Angers, 1979, Sijthoff \&\
  Noordhoff, Alphen aan den Rijn, 1980, pp.~273--310.

\bibitem[SS]{SS}
S.~P. Smith and J.~T. Stafford, \emph{Regularity of the four-dimensional
  {S}klyanin algebra}, Compositio Math. \textbf{83} (1992), no.~3, 259--289.

\bibitem[St1]{St1}
Darin~R. Stephenson, \emph{Artin-{S}chelter regular algebras of global
  dimension three}, J. Algebra \textbf{183} (1996), no.~1, 55--73.

\bibitem[St2]{St2}
\bysame, \emph{Algebras associated to elliptic curves}, Trans. Amer. Math. Soc.
  \textbf{349} (1997), no.~6, 2317--2340.

\bibitem[St3]{Ste}
\bysame, \emph{The geometry of noncommutative graded algebras}, preliminary
  version, 1998.
  
\bibitem[SZ]{SZ}
Darin~R. Stephenson and James~J. Zhang, \emph{Growth of graded {N}oetherian
  rings}, Proc. Amer. Math. Soc. \textbf{125} (1997), no.~6, 1593--1605.
  
\bibitem[V]{V}
James~S. Vandergraft, \emph{Spectral properties of matrices which have
  invariant cones}, SIAM J. Appl. Math. \textbf{16} (1968), 1208--1222.

\bibitem[W]{W}
Joachim Wehler, \emph{${K}3$-surfaces with {P}icard number $2$}, Arch. Math.
  (Basel) \textbf{50} (1988), no.~1, 73--82.

\end{thebibliography}

\providecommand{\bysame}{\leavevmode\hbox to3em{\hrulefill}\thinspace}


\end{document}